\documentclass{article}

\title{Using the VBARMS method in parallel computing}

\author{
Bruno Carpentieri\footnote{Institute of Mathematics and Computing Science - University of Groningen, 9747 AG Groningen, The Netherlands - e-mail: \url{b.carpentieri@rug.nl},
\url{j.liao@rug.nl}}
    \and  
Jia Liao$^\fnsymbol{footnote}$
    \and 
Masha Sosonkina\footnote{Department of Modeling, Simulation \& Visualization Engineering - Old Dominion University,  Norfolk, VA 23529 - e-mail: \url{msosonki@odu.edu}}
	\and 
Aldo Bonfiglioli\footnote{Scuola di Ingegneria - University of Basilicata, Potenza, Italy - e-mail:
\url{aldo.bonfiglioli@unibas.it}}	
}

\date{}

\usepackage{hyperref}
\usepackage{graphicx,epsfig}
\usepackage{amsmath,amssymb}
\usepackage{placeins}
\usepackage[noend]{algorithmic}
\usepackage{algorithm}
\usepackage[scriptsize]{subfigure}
\usepackage{url}
\usepackage{cite}
\usepackage{listings}
\usepackage{verbatim}
\usepackage{colortbl}
\usepackage{subfigmat}

\usepackage{upgreek}
\usepackage{longtable}
\usepackage{multirow}

\newcommand{\ms}{\medskip}
\newcommand{\bs}{\bigskip}
\newcommand{\tabincell}[2]{\begin{tabular}{@{}#1@{}}#2\end{tabular}}

\newcommand{\R}{\mathbb{R}}

\newcommand{\N}{\mathbb{N}}

\begin{document}

\bibliographystyle{plain}

\maketitle

\makeatother

\renewcommand{\thefootnote}{\arabic{footnote}}
\newcommand{\eps}{\varepsilon}
\newcommand{\bb}{\phantom{8}}
\newsavebox{\fmbox}
\newenvironment{fmpage}[1]
{ \begin{lrbox}
    {\fmbox}
    \begin{minipage}{#1}}
    {\end{minipage}
\end{lrbox}
\fbox{\usebox{\fmbox}} }

\newenvironment{maliste}%
{ \begin{list}%
    {$\bullet$}%
    {\setlength{\labelwidth}{30pt}%
         \setlength{\leftmargin}{35pt}%
     \setlength{\itemsep}{\parsep}}}%
{ \end{list} }

\begin{abstract}

The paper describes an improved parallel MPI-based implementation of VBARMS, a variable block variant of the pARMS preconditioner proposed by Li,~Saad and Sosonkina [NLAA, 2003] for solving general nonsymmetric linear systems. The parallel VBARMS solver can detect automatically exact or approximate dense structures in the linear system, and exploits this information to achieve improved reliability and increased throughput during the factorization. A novel graph compression algorithm is discussed that finds these approximate dense blocks structures and requires only one simple to use parameter. A complete study of the numerical and parallel performance of parallel VBARMS is presented for the analysis of large turbulent Navier-Stokes equations on a suite of three-dimensional test cases.

\medskip 

\end{abstract}

\section{Introduction}~\label{sec:1}

The initial motivation for this study is the design of robust preconditioning techniques for solving sparse block structured linear systems arising from the finite element / finite volume analysis of turbulent flows in computational fluid dynamics applications. Over the last few years we have developed block multilevel incomplete LU (ILU) factorization methods for this problem class, and we have found them very effective in reducing the number of GMRES iterations compared to their pointwise analogues~\cite{VBARMS}. This class of preconditioners can offer higher parallelism and robustness than standard ILU algorithms especially for solving large problems, thanks to their multilevel  mechanism. Exploiting existing block structures in the matrix can help reduce numerical instabilities during the factorization and achieve higher flops to memory ratios on modern cache-based computer architectures. Sparse matrices arising from the solution of systems of partial differential equations often exhibit perfect block structures consisting of fully dense (typically small) nonzero blocks in their sparsity pattern, e.g., when several unknown physical quantities are associated with the same grid point. For example, a plane elasticity problem has both $x$- and $y$-displacements at each grid point; a Navier-Stokes system for turbulent compressible flows would have five distinct variables (the density, the scaled energy, two components of the scaled velocity, and the turbulence transport variable) assigned to each node of the physical mesh; a bidomain system in cardiac electrical dynamics couples the intra-and extra-cellular electric potential at each ventricular cell of the heart. Upon numbering consecutively the $\ell$ distinct variables associated with the same grid point, the permuted matrix has a sparse block structure with nonzero blocks of size $\ell \times \ell$. The blocks are fully dense if variables at the same node are mutually coupled. 

Our recently developed variable block algebraic recursive multilevel solver (shortly, VBARMS) can detect fine-grained dense structures in the linear system automatically, without any user's knowledge of the underlying problem, and exploit them efficiently during the factorization~\cite{VBARMS}. Preliminary experiments with a parallel MPI-based implementation of VBARMS for distributed memory computers, presented in a conference contribution~\cite{ppam2014}, showed the robustness of the proposed method for solving some larger matrix problems arising in different fields. In this paper, capitalizing on those results, we introduce a new graph-based compression algorithm to construct the block ordering in VBARMS, which extends the method proposed by Ashcraft in~\cite{ashc:95} and requires only one simple to use parameter (Section~\ref{sec:2}); we describe in Section~\ref{sec:3} a novel implementation of the block partial factorization step that proves to be noticeably faster than the original one presented in~\cite{VBARMS}; finally, in Section~\ref{sec:4}, we assess the parallel performance of our parallel VBARMS code for solving turbulent Navier-Stokes equations in fully coupled form on large realistic three-dimensional meshes; in the new parallel implementation, we use a parallel graph partitioner to reduce the graph partitioning time  significantly compared to the experiments presented in~\cite{ppam2014}.

\section{Graph compression techniques}~\label{sec:2}

It is known that block iterative methods often show faster convergence rate than their pointwise analogues in the solution of many classes of two- and three-dimensional partial differential equations (PDEs). When the domain is discretized by cartesian grids, a regular partition may also provide an effective matrix partitioning. For example, in the case of the simple Poisson's equation with Dirichlet boundary conditions, defined on a rectangle $(0,\ell_1)\times(0,\ell_2)$ discretized uniformly by $n_1+2$ points in the interval $(0,\ell_1)$ and $n_2+2$ points in $(0,\ell_2)$, upon numbering the interior points in the natural ordering by lines from the bottom up, one obtains a $n_2\times n_2$ block tridiagonal matrix with square blocks of size $n_1\times n_1$; the diagonal blocks are tridiagonal matrices and the off-diagonal blocks are diagonal matrices. For large finite element discretizations, it is common to use substructuring, where each substructure of the physical mesh corresponds to one sparse block of the system. However, if the domain is highly irregular or the matrix does not correspond to a differential equation, finding the best block partitioning is much less obvious. In this case, graph reordering techniques are worth considering.

The PArameterized BLock Ordering (PABLO) method proposed by O'Neil and Szyld is one of the first matrix partitioning algorithms specifically designed for block iterative solvers~\cite{nesz:90}.  
The algorithm selects groups of nodes in the adjacency graph of the coefficient matrix such that the corresponding diagonal blocks are either full or very dense. It has been shown that classical block stationary iterative methods such as block Gauss-Seidel and SOR methods combined with the PABLO ordering require fewer operations than their point analogues for the the finite element discretization of a Dirichlet problem on a graded L-shaped region, as well as on the 9-point discretization of the Laplacian operator on a square grid. The complexity of the PABLO algorithm is proportional to the number of nodes and edges in both time and space. 

Another useful approach for blocking a matrix $A$ is to find block independent sets in the adjacency graph of $A$~\cite{ARMS}. A block independent set is defined as a set of groups of nodes (or unknowns) having the property that there is no coupling between nodes of any two different groups, while nodes within the same group may be coupled. Independent sets of unknowns in a linear system can be eliminated simultaneously at a given stage of Gaussian Elimination. For this reason, this type of oredering is extensively adopted in linear solvers design. Independent sets may be computed by using simple graph algorithms which traverse the vertices of the adjacency graph of $A$ in the natural order $1,2,\ldots,n$, mark each visited vertex $v$ and all of its nearest neighbors connected to $v$ by an edge, and add $v$ and each visited node that is not already marked to the current independent set partition~\cite{Saad:1996:IME}. Upon renumbering nodes one partition after the other, followed as last by interface nodes straddling between separate partitions, one obtain a permutation of $A$ in the form
\begin{equation}~\label{eq:indset1}
P A P^T  = \left( {\begin{array}{*{20}c}
   D & F  \\
   E & C  \\
\end{array}} \right),
\end{equation}
where $D$ is a block diagonal matrix. The nested dissection ordering by George~\cite{geli:81}, mesh partitioning, or further information from the set of nested finite element grids of the underlying problem~\cite{OAxelsson_PSVassilevski_1989a,OAxelsson_PSVassilevski_1990a,NGILU} can be used as an alternative to the greedy independent set algorithm described above. Additionally, the numerical values of $A$ may be incorporated in the ordering to produce more robust factorizations~\cite{ARMS}. 

However, finite element and finite difference matrices often possess also fine-grained block structures that can be exploited in iterative solvers. If there is more than one solution component at a grid point, the corresponding matrix entries may form a small dense block and optimized codes can be used for dense factorizations in the construction of the preconditioner and dense matrix-vector products in the sparse matrix-vector product operation for  better performance, see e.g.~\cite{VBARMS,IMF,guge:10,VBILU,vuduc:2007}. A block incomplete LU factorization (ILU) method is one preconditioning technique that treats small dense submatrices of $A$ as single entities, and the VBARMS method discussed in this paper can be seen as its natural multilevel generalization. An important advantage of block ILU versus conventional ILU is the potential gain obtained from using optimized level 3 basic linear algebra subroutines (BLAS3). Column indices and pointers can be saved by storing the matrix as a collection of blocks using the variable block compressed sparse row (VBCSR) format, where each value in the CSR format is a dense array. On indefinite problems, computing with blocks instead of single elements enables us a better control of pivot breakdowns, near singularities, and other sources of numerical instabilities. These facts have been assessed in our previous contribution~\cite{VBARMS}.

The method proposed by Ashcraft in~\cite{ashc:95} is one of the first compression techniques for finding dense blocks in the sparsity pattern of a matrix. The algorithm searches for sets of rows or columns having the exact same pattern. From a graph viewpoint, it looks for vertices of the adjacency graph 
$(V,E)$ of $A$ having the same adjacency list. These are also called {\it{indistinguishable nodes}} or {\it{cliques}}.
The algorithm assigns a {\it{checksum}} quantity to each vertex, e.g., using the function
\begin{equation}\label{checksum}
chk(u) = \sum\limits_{(u,w) \in E} w ,
\end{equation}
and then sorts the vertices by their checksums. This operation takes $|E| + |V|\log |V|$ time. If $u$ and $v$ are indistinguishable, then $chk(u) = chk(v)$. Therefore, the algorithm examine nodes having the same checksum to see if they are indistinguishable. The ideal checksum function would assign a different value for each different row pattern that occurs but it is not practical because it may quickly lead to huge numbers that may not even be machine-representable. Since the time cost required by Ashcraft's method is generally negligible relative to the time it takes to solve the system, simple checksum functions such as~\eqref{checksum} are used in practice~\cite{ashc:95}. 

Sparse unstructured matrices may sometimes exhibit {\it{approximate dense blocks}} consisting mostly of nonzero entries except only a few zeros inside the blocks. By treating these few zeros as nonzero elements, with a little sacrifice of memory, a block ordering may be generated for an iterative solver. 
Computing approximate dense structures may enable us to enlarge existing blocks and to use BLAS3 operations more efficiently in the iterative solution, but it may also increase the memory costs and the probability to encounter singular blocks during the factorization~\cite{VBARMS}. Two important performance measures to gauge the quality of the block ordering computed are the {\it{average block density}} (${av\_bd}$) value, defined as the amount of nonzeros in the matrix divided by the amount of elements in the nonzero blocks, and the {\it{average block size}} (${av\_bs}$) value, which is the ratio between the sum of dimensions of the square diagonal blocks divided by the number of diagonal blocks. From our computational experience, high average block density values around 90\% are necessary to prevent the occurrence of singular blocks during the factorization. 

\ms 

\subsection{The angle-based method}

Approximate dense blocks in a matrix may be computed by numbering consecutively rows and columns having a similar nonzero structure. However, this would require a new checksum function that preserves the proximity of patterns, in the sense that close patterns would result in close checksum values. Unfortunately, this property does not hold true for Ashcraft's algorithm in its original form. 
In~\cite{VBILU}, Saad proposed to compare angles of rows (or columns) to compute approximate dense structures in a matrix $A$. Let $C$ be the pattern matrix of $A$, which by definition has the same pattern as $A$ and nonzero values equal to one. The method proposed by Saad computes the upper triangular part of $C C^T$. Entry $(i,j)$ is the inner product (the cosine value) between row $i$ and row $j$ of $C$ for $j>i$. 
A parameter $\tau$ is used to gauge the proximity of row patterns. If the cosine of the angle between rows $i$ and $j$ is smaller than $\tau$, row $j$ is added to the group of row $i$. For $\tau=1$ the method will compute perfectly dense blocks, while for $\tau<1$ it may compute larger blocks where some zero entries are padded in the pattern. To speed up the search, it may be convenient to run a first pass with the checksum algorithm to detect rows having an identical pattern, and group them together; then, in a second pass, each non-assigned row is scanned again to determine whether it can be added to an existing group. In practice, however, it may be difficult to predict the average block density obtained using a given value of $\tau$ . For example, the experiments reported in Table~\ref{tab:bdsaad} show that $\tau=0.58$ returns a block density of $86.37\%$ for the {VENKAT01} matrix and of $45.06\%$ for the {STACOM} matrix. 

\begin{table}[!ht]
\centering
\begin{tabular}{p{50px}|c|c|c|c|c}
Matrix   & $\tau=0.56$ & $\tau=0.57$ & $\tau=0.58$ & $\tau=0.59$ & $\tau=0.60$\\ \hline
{STACOM}   & 25.63 & 25.68 & 45.06 & 50.83 & 52.02\\ \hline
{K3PLATES} & 37.78 & 38.73 & 58.62 & 58.70 & 59.16\\ \hline 
{OILPAN}   & 50.08 & 50.09 & 50.23 & 50.23 & 90.65\\ \hline
{VENKAT01} & 29.71 & 29.71 & 86.37 & 86.37 & 86.37\\ \hline
{RAE}      & 26.40 & 26.48 & 49.48 & 50.71 & 51.96
\end{tabular}\\[\baselineskip]
\begin{tabular}{p{50px}|c|c|c|c|c}
Matrix   & $\tau=0.64$ & $\tau=0.65$ & $\tau=0.66$ & $\tau=0.67$ & $\tau=0.68$\\ \hline
{RAEFSKY3} & 63.32 & 63.32 & 63.32 & 95.23 & 95.23\\ \hline
{BMW7ST\_1} & 49.29 & 50.11 & 50.66 & 68.85 & 74.00\\ \hline
{S3DKQ4M2} & 64.29 & 64.29 & 64.29 & 97.52 & 97.52\\ \hline
{PWTK}     & 57.05 & 57.31 & 57.48 & 94.23 & 94.75
\end{tabular}
\caption{Average block density value (\%) obtained from the angle compression algorithm for different values of $\tau$. \label{tab:bdsaad}}
\end{table}

The cost of Saad's method is closer to that of checksum-based methods for cases in which a good blocking already exists, and in most cases it remains inferior to the cost of the least expensive block LU factorization, i.e., block ILU(0).

\subsection{Graph-based compression}

We revisited Saad's angle-based method to develop a new compression algorithm that computes a block ordering having an average block density ${av\_bd}$ not smaller than a user-specified value $\mu$. This may simplify the parameter selection procedure. 
The method proceeds in two steps. First, using the checksum algorithm it groups rows having equal nonzero structure and builds the quotient graph $G/{\mathcal{B}} = (V_\mathcal{B}, E_\mathcal{B})$. In  $G/{\mathcal{B}}$, nodes corresponding to rows with identical pattern are coalesced into one single node of $V_\mathcal{B}$ (also called supernode or supervertex). An edge connects supervertices $Y$ and $Z$ of $V_\mathcal{B}$ if there exists an edge in $G = (V, E)$ connecting a vertex in $Y$ to a vertex in $Z$. If $A$ is unsymmetric, we assume to operate on the symmetrized graph of $A+A^T$; thus the edge orientation is not important. Afterwards, the algorithm merges pairs of supernodes $(Y,Z)$, for $Z$ adjacent to $Y$ in $G/{\mathcal{B}}$, provided that the average block density value $av\_bd$ of the new block ordering after this operation does not drop below $\mu$. Otherwise, the algorithm will stop to prevent near-singularities during the block factorization. The total size of the rows and columns spanned by this new block is
\[
T = 2 \cdot \left| {adj(Y) \cup adj(X)} \right| \cdot \left| {Y \cup X} \right| - \left| {Y \cup X} \right|^2 ,
\]
which is the amount of nonzero rows and columns times the size of the supernode minus the square block on the diagonal which we count twice since we count both columns and rows. The nonzeros spanned by the new block is
\[
N = 2 \cdot \sum\limits_{Z \in Y \cup X} {\left| {adj(Z)} \right|}  - \sum\limits_{Z \in Y \cup X} {\left| {adj(Z) \cap \left( {Y \cup X} \right)} \right|} ,
\]
which is the amount of adjacent nodes per node inside the supernode minus the amount of nodes inside the diagonal block, which is again counted twice. The complete graph-based algorithm is sketched in Algorithm~\ref{fig:graphcompr}. It requires only one simple to use parameter $\mu$. If we desire a block ordering having an average block density around $60\%$, we simply set $\mu=0.6$. In contrast, a correct tuning of $\tau$ may require to run the full solver to see if a singular block is encountered during the factorization.

\subsection{Experiments}

In Table~\ref{tab:compr3} we give some comparative performance figure to show the viability of the graph algorithm. In our runs we attempted to find the optimal value of $\tau$ by trial and error. By optimal value we mean the one that minimizes the number of GMRES iterations required to reduce the initial residual by $6$ orders of magnitude using a standard block incomplete LU factorization as a preconditioner for GMRES. The optimal value for the parameter $\tau$ was calculated by running the angle algorithm with different $\tau \in [0.5,1.0]$, by increments of $0.1$ at every run. The results evidence the difficulty to compute a unique value which is nearly optimal for every problem. On the other hand, for the graph method we set $\mu=0.7$ which gave us a minimum block density of $70\%$ for every matrix. We see that the new compression algorithm is very competitive and additionally may be simple to use. In Tables~\ref{tab:compr3} we also report on the timing to compute the block ordering by both compression techniques, and for solving the linear system. The new graph algorithm is in most cases up to three times slower than the angle algorithm. However, this is not a big downside because the compression time is considerably smaller than the total solution time, and computing the optimal value of $\tau$ may require several runs as we explained. Clearly, the compression time increases when $\mu$ decreases since we merge more supernodes in this circumstance. By the way, both compression methods helped reduce iterations. Without blocking, no convergence was achieved in 1000 iterations using pointwise ILUT on the OILPAN, K3PLATES, S3DKQ4M2, OLAFU, RAE, NASASRB, CT20STIF, RAEFSKY3, BCSSTK35, STACOM problems at equal or higher memory usage. On the other hand, no evident gain was observed from using level-2 BLAS routines in the sparse matrix-vector product operation, probably due to the small block size.

\begin{table}[!ht]
 \begin{center}
  \begin{tabular}{l|r|l|r|c}
    Name & Size & Application & nnz(A) & symmetry \\
    \hline
	OILPAN & 73752 & Structural problem & 2148558 & symmetric value \\
	K3PLATES & 11107 & FE stiffness matrix & 378927 & symmetric value \\
	VENKAT01 & 62424 & Unstructured 2D Euler solver & 1717792 & symmetric structure\\
	PWTK & 217918 & Pressurized wind tunnel & 11524432 & symmetric value \\
	S3DKQ4M2 & 90449 & Structural mechanics & 2455670 & symmetric value \\
	OLAFU & 16146 & Structural problem & 1015156 & symmetric value \\
	RAE & 52995 & Turbulence analysis & 1748266 & symmetric structure \\
	BMW7ST\_1 & 141347 & Stiffness matrix & 7318399 & symmetric value \\
	NASASRB & 54870 & Shuttle rocket booster & 2677324 & symmetric value  \\
	CT20STIF & 52329 & Stiffness matrix engine block & 2600295 & symmetric value \\
	RAEFSKY3 & 21200 & Fluid structure interaction turbulence problem & 1488768 & symmetric structure \\
	HEART1 & 3557 & Quasi-static FEM of a heart & 1385317 & symmetric structure \\
	BCSSTK35 & 30237 & Automobile seat frame & 1450163 & symmetric value \\
	STACOM & 8415 & Compressible flow & 271936 & symmetric structure \\
  \end{tabular}
 \caption{\label{tab:matrices}Set and characteristics of test matrix problems.}
 \end{center}
\end{table}

\begin{longtable}{p{50px}|p{38px}|c|c|c|c|c|p{38px}|p{25px}}
\hline
{Matrix} & {Method} & {$\tau/ \mu$} & {$av\_bd$~(\%)}  & {$av\_bs$} & \begin{tabular}{c}Blocking \\ time (s)\end{tabular} & \begin{tabular}{c}Solving \\ time (s)\end{tabular} & Mem & Its\\
\hline
\multirow{2}{*}{OILPAN} & Angle & 0.70 & 95.94 & 7.36 & 0.03 & 4.18 & 0.26 & 198\\ 
 & Graph & 0.70 & 95.02 & 7.42 & 0.08 & 4.17 & 0.27 & 198\\ 
\hline
\multirow{2}{*}{K3PLATES} & Angle & 0.60 & 59.16 & 7.90 & 0.00 & 0.7 & 0.3 & 239\\ 
 & Graph & 0.70 & 89.50 & 5.65 & 0.01 & 0.7 & 0.18 & 241\\ 
\hline
\multirow{2}{*}{VENKAT01} & Angle & 0.70 & 99.94 & 4.00 & 0.02 & 0.43 & 1.33 & 9\\ 
 & Graph & 0.70 & 94.05 & 4.28 & 0.08 & 0.48 & 1.58 & 9\\ 
\hline
\multirow{2}{*}{PWTK} & Angle & 0.60 & 56.95 & 12.17 & 0.09 & 26.38 & 6.85 & 117\\ 
 & Graph & 0.70 & 78.16 & 7.31 & 0.35 & 32.64 & 4.5 & 137\\ 
\hline
\multirow{2}{*}{S3DKQ4M2} & Angle & 1.00 & 100.00 & 5.93 & 0.03 & 9.57 & 1.09 & 214\\ 
 & Graph & 0.70 & 77.92 & 7.81 & 0.12 & 15.1 & 1.42 & 309\\ 
\hline
\multirow{2}{*}{OLAFU} & Angle & 0.80 & 81.75 & 6.47 & 0.02 & 1.2 & 3.14 & 54\\ 
 & Graph & 0.70 & 79.66 & 6.58 & 0.11 & 1.63 & 3.75 & 57\\ 
\hline
\multirow{2}{*}{RAE} & Angle & 0.80 & 95.83 & 4.67 & 0.03 & 8.85 & 9.53 & 49\\ 
 & Graph & 0.70 & 86.21 & 4.64 & 0.13 & 15.74 & 13.8 & 42\\ 
\hline
\multirow{2}{*}{BMW7ST\_1} & Angle & 0.70 & 77.16 & 7.28 & 0.08 & 0.35 & 0.18 & 5\\ 
 & Graph & 0.70 & 79.54 & 6.65 & 0.29 & 0.48 & 0.17 & 9\\ 
\hline
\multirow{2}{*}{NASASRB} & Angle & 0.80 & 90.87 & 4.24 & 0.05 & 7.51 & 5.23 & 30\\ 
 & Graph & 0.70 & 77.62 & 4.20 & 0.20 & 12.39 & 7.46 & 16\\ 
\hline
\multirow{2}{*}{CT20STIF} & Angle & 0.70 & 66.05 & 6.55 & 0.04 & 0.69 & 0.18 & 44\\ 
 & Graph & 0.70 & 78.42 & 4.76 & 0.16 & 1.18 & 0.14 & 56\\ 
\hline
\multirow{2}{*}{RAEFSKY3} & Angle & 0.70 & 95.23 & 8.63 & 0.01 & 0.08 & 0.13 & 13\\ 
 & Graph & 0.70 & 77.67 & 10.56 & 0.02 & 0.09 & 0.17 & 15\\ 
\hline
\multirow{2}{*}{HEART1} & Angle & 0.90 & 98.81 & 18.62 & 0.00 & 0.5 & 0.78 & 151\\ 
 & Graph & 0.70 & 0.00 & 0.00 & 0.00 & - & - & -\\ 
\hline
\multirow{2}{*}{BCSSTK35} & Angle & 0.60 & 51.95 & 11.03 & 0.01 & 2.1 & 0.29 & 209\\ 
 & Graph & 0.70 & 78.72 & 6.57 & 0.05 & 2.66 & 0.18 & 235\\ 
\hline
\multirow{2}{*}{STACOM} & Angle & 0.90 & 97.00 & 4.36 & 0.00 & 0.25 & 5.19 & 31\\
 & Graph & 0.70 & 84.51 & 4.47 & 0.01 & 0.29 & 5.65 & 33\\ \hline 
\caption{Experiments with the angle-based and the graph-based compression methods. The optimal value of $\tau$ is used for the angle-based algorithm. The value $\mu=0.7$ is used for the graph-based algorithm in all our runs.}\label{tab:compr3}
\end{longtable}

\bs 

\bs 

\begin{algorithm}
\caption{\label{fig:graphcompr}{\it The graph based compression algorithm}.}
\begin{algorithmic}[1]
\STATE Compute the keys $k_i=chk(i)$ for all vertices $i \in V =\{1,\ldots,n\}$
\STATE Set processed nodes $p_i=0 ~ \forall i=1,\ldots,n$
\STATE Make a set of supernodes $\mathcal{V}=\emptyset$
\STATE Set $s$ to the indices $V$ sorted by the corresponding value in $k$
\FOR{$i=s_1,\ldots,s_n$}
  \IF{$p_i \ne 1$} 
  \STATE Add a new supernode $Y_i$ to $\mathcal{V}$
  \FOR{$j=s_{i+1},\ldots,s_n$}
	  \IF{$k_i \ne k_j$} 
		  \STATE{\bf{break}}
	  \ENDIF
	  \IF{$adj(i) = adj(j)$} 
		  \STATE Add node $j$ to $Y_i$
		  \STATE Set $p_j=1$
	  \ENDIF
  \ENDFOR
  \ENDIF
\ENDFOR
\STATE Make a map $\mathcal{M}:i \mapsto \left\{ {Z \in \mathcal{V} | ~ i \in adj(Z)} \right\}$
\FOR{$X \in \mathcal{V}$}
	\FOR{$Z \in \bigcup\nolimits_{i \in X} {\mathcal{M}(i)} $}
	    \STATE Update the average block density value $av\_bd$ after merging $X$ and $Z$
		\IF{$av\_bd \ge \mu$}
		  \STATE{$X = X \cup Z$}
		  \STATE{$\mathcal{V} = \mathcal{V}  \backslash Z$}
		\ENDIF
	\ENDFOR
\ENDFOR
\end{algorithmic}
\end{algorithm}

\medskip 

\begin{figure}[!ht]
 \centering
 \subfigure[Using the PABLO algorithm]
 { \label{fig:ex1}
   \includegraphics[width=7cm]{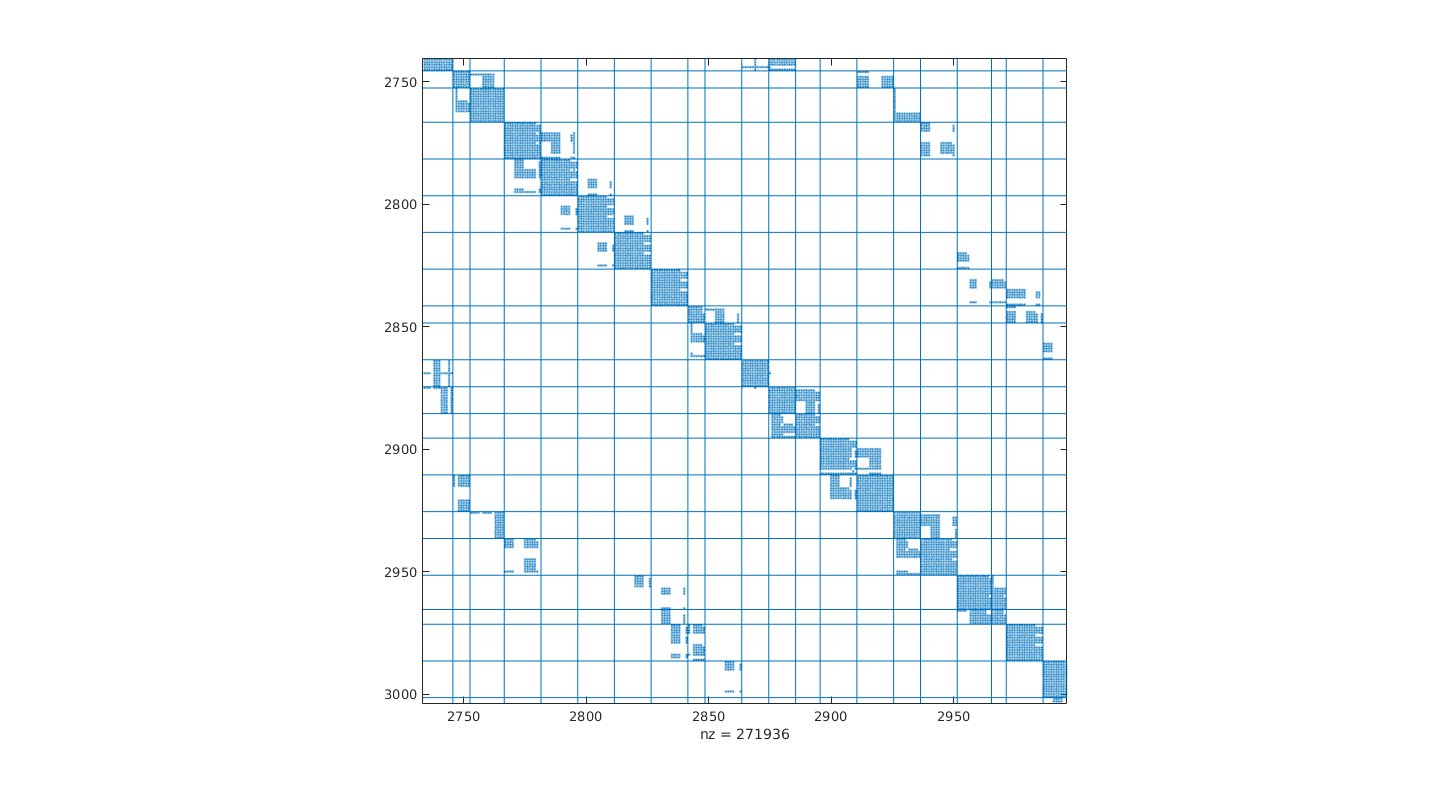}
 }
 \subfigure[Using the graph algorithm]
 { \label{fig:ex2}
   \includegraphics[width=7cm]{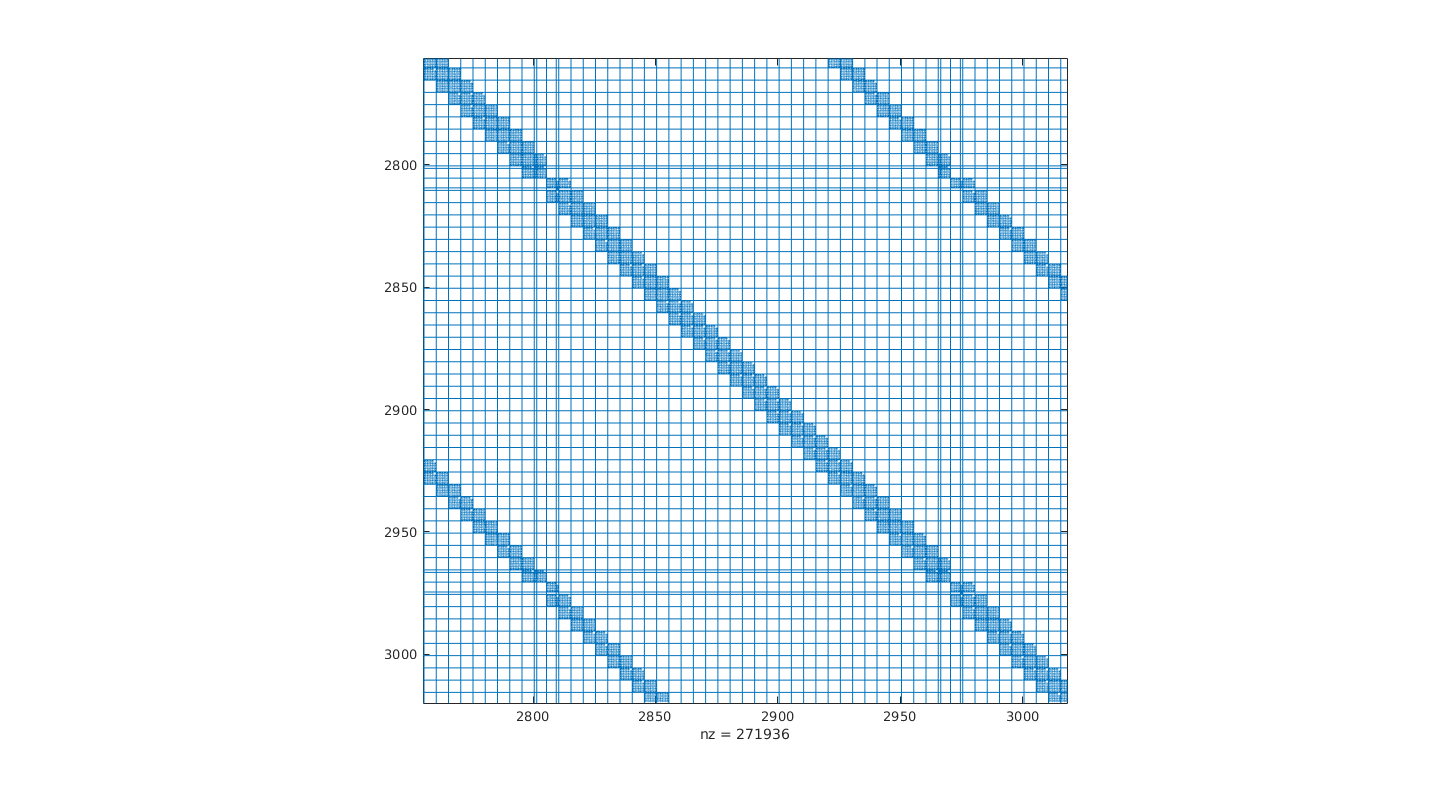}
 }
 \caption{\label{fig:PABLO_PATTERN} Block patterns computed by different compression methods for the STACOM problem.}
\end{figure}

For the sake of comparison, we also ran some experiments using the PABLO algorithm introduced by O'Neil and Szyld in~\cite{nesz:90}, in combination with block incomplete LU factorization preconditioning. The convergence results are reported in Table~\ref{tab:PABLO_VBILUT}, and a comparison of patterns produced by the two compression techniques is shown in Figure~\ref{fig:PABLO_PATTERN} for two matrices. We observe that the block ordering computed by PABLO may produce larger blocks compared to the graph and angle methods. However, the average block size can be significantly smaller, probably due to the design philosophy of PABLO that attempts to maximize the density of the diagonal blocks of a matrix. The convergence results show that overall the resulting block ordering may be less suitable for block factorization.

\begin{table}[h]
\begin{center}
\begin{tabular}{c|c|c|c|c|c}
\hline
Matrix  & $av\_bd$ & $av\_bs$ & \begin{tabular}{c}Total \\ time (s)\end{tabular} & Mem & Its\\
\hline
STACOM  & 66.54 & 2.38 & 6.22 & 11.02 & 152\\
\hline
K3PLATES  & 83.51 & 2.00 &  8.94 & 5.54& 329\\
\hline
OLAFU  & 89.60 & 2.00 & 7.66 & 3.89 & 84\\
\hline
RAE  & 68.28 & 2.34 & 412.89  & 26.75 & 1000\\
\hline
\end{tabular}
\end{center}
\caption{Pablo performance and Pablo with VBILUT. Block density refers to the average block density of the block ordering, Block size is the average block size, Total time includes the preconditioning construction and the solving time, Mem is the ratio between the number of nonzeros in the preconditioner and in the matrix.}\label{tab:PABLO_VBILUT}
\end{table}

\section{The VBARMS method}~\label{sec:3}

The VBARMS method discussed in this paper incorporates compression techniques to maximize computational efficiency during the factorization. We recall briefly below the main steps of the algorithm and we point the reader to~\cite{VBARMS} for further details. After permuting the coefficient matrix $A$ in block form as
\begin{equation} \label{eq:blockperm}
\widetilde{A} \approx P_B AP_B^T = \left[ {\begin{array}{*{20}c}
   {\widetilde A_{11} } & {\widetilde A_{12} } &  \cdots  & {\widetilde A_{1p} }  \\
   {\widetilde A_{21} } & {\widetilde A_{22} } &  \cdots  & {\widetilde A_{2p} }  \\
    \vdots  &  \vdots  &  \ddots  &  \vdots   \\
   {\widetilde A_{p1} } & {\widetilde A_{p2} } &  \cdots  & {\widetilde A_{pp} }  \\
\end{array}} \right],
\end{equation}
where the diagonal blocks $\widetilde A_{ii}$, $i=1,\ldots,p$ are $n_i \times n_i$, the off-diagonal blocks $\widetilde A_{ij}$ are  $n_i \times n_j$, and $P_B$ is the permutation matrix of the block ordering computed by the compression algorithm,
we can represent the adjacency graph of $\widetilde{A}$ by the quotient graph of ${A}+{A^T}$~\cite{geli:81}, which is smaller. Let $\mathcal{B}$ the partition into blocks given by~(\ref{eq:blockperm}). The quotient graph $\mathcal{G}/\mathcal{B}=(V_\mathcal{B},E_\mathcal{B})$ is constructed by coalescing the vertices of each block $\widetilde A_{ii}$, for $i=1,\ldots,p$, into one supervertex (or supernode) $Y_i$. We can write
\[
V_\mathcal{B}  = \left\{ {Y_1 , \ldots ,Y_p } \right\}, ~~~ E_\mathcal{B}  = \left\{ {\left( {Y_i ,Y_j } \right)~|~\exists v \in Y_i ,w \in Y_j~\text{s.t.}~(v,w) \in E} \right\}.
\]
where $(V,E)$ is the graph of $A+A^T$. An edge connects two supervertices $Y_i$ and $Y_j$ if there exists an edge of $(V,E)$ connecting a vertex of the block $A_{ii}$ to a vertex of the block $A_{jj}$. 

\ms

The complete pre-processing and factorization process of VBARMS consists of the following steps.

\begin{description}
\item [Step~1] Using the angle-based or the graph-based compression algorithms described in Section~\ref{sec:2}, compute a block ordering $P_B$ of $A$ such that, after permutation, the matrix $P_B A P_B^T$ has fairly dense nonzero blocks. 

\item [Step~2] Scale the matrix permuted at {\bf{Step~1}} as  $S_1 P_B A P_B^T S_2$, where $S_1$ and $S_2$ are two diagonal matrices such that the 1-norm of the largest entry in each row and column becomes smaller or equal than one.

\item [Step~3] Apply the block independent sets (or the nested dissection) algorithms to the quotient graph $\mathcal{G}/\mathcal{B}$ and compute an independet sets ordering $P_{I}$ of $\mathcal{G}/\mathcal{B}$. 
Upon permutation by $P_{I}$, the matrix obtained at {\bf{Step~2}} will write as
\begin{equation}~\label{eq:indset1}
P_{I} S_1 P_B A P_B^T S_2 P_{I}^T  = \left( {\begin{array}{*{20}c}
   D & F  \\
   E & C  \\
\end{array}} \right).
\end{equation}
We use a simple weighted greedy algorithm for computing the ordering $P_{I}$~\cite{ARMS}.

In the $2 \times 2$ partitioning~(\ref{eq:indset1}), the upper left-most matrix $D \in \R^{m \times m}$ is  block diagonal like in ARMS. However, due to the block permutation ({\bf{Step~1}}), the diagonal blocks $D_i$ of $D$ are block sparse matrices while in ARMS they are sparse unstructured. The matrices $F \in \R^{m \times (n-m)}$, $E \in \R^{(n-m) \times m}$, 
$C \in \R^{(n-m) \times (n-m)}$ are also block sparse, because of the same reason. 

\item [Step~4] Factorize the matrix in~(\ref{eq:indset1}) as
\begin{equation}\label{eq:vbarms}
\left( {\begin{array}{*{20}c}
    D &  F  \\
    E &  C  \\
\end{array}} \right) = \left( {\begin{array}{*{20}c}
    L & 0  \\
   { E  U^{ - 1} } & I  \\
\end{array}} \right) \times \left( {\begin{array}{*{20}c}
    U & { L^{ - 1}  F}  \\
    0 & { A_1 }  \\
\end{array}} \right),
\end{equation}
where $I$ is the identity matrix of appropriate size, and 
\begin{equation}\label{eq:schurcomp}
	A_1 = C- E  D^{-1}  F.
\end{equation}
is the Schur complement corresponding to $C$. Observe that the Schur complement is also block sparse and it has the same block structure as matrix $C$.

\end{description}

{\bf{Steps~2-4}} can be repeated on the reduced system a few times until the Schur complement is small enough. Denoting by $A_\ell$ the reduced Schur complement matrix at level $\ell$, for $\ell > 1$, after scaling and preordering $A_\ell$ a system with coefficient matrix
\begin{equation}\label{eq:sysl}
	P_I^{(\ell)}  D_1^{(\ell)} A_\ell  D_2^{(\ell)} (P_I^{(\ell)})^T = \left( {\begin{array}{*{20}c}
    D_\ell &  F_\ell  \\
    E_\ell &  C_\ell  \\
\end{array}} \right)
 = \left( {\begin{array}{*{20}c}
    L_\ell & 0  \\
   { E_\ell  U_\ell^{ - 1} } & I  \\
\end{array}} \right) \times \left( {\begin{array}{*{20}c}
    U_\ell & { L_\ell^{ - 1}  F_\ell}  \\
    0 & { A_{\ell+1} }  \\
\end{array}} \right)
\end{equation}
needs to be solved, with
$D_\ell \in \R^{m_\ell \times m_\ell}$, $F_\ell \in \R^{m_\ell \times (n_\ell-m_\ell)}$, $E_\ell \in \R^{(n_\ell-m_\ell) \times m_\ell}$, $C_\ell \in \R^{(n_\ell-m_\ell) \times (n_\ell-m_\ell)}$, and 
\begin{equation}\label{eq:schurl}
A_{\ell+1} = C_\ell - E_\ell  D_\ell^{-1}  F_\ell \in \R^{(n_\ell-m_\ell) \times (n_\ell-m_\ell)}.
\end{equation}
Calling  
\[
x_\ell = \left( {\begin{array}{c}
    y_\ell   \\
    z_\ell   \\
\end{array}} \right), ~~~ b_\ell = \left( {\begin{array}{c}
    f_\ell   \\
    g_\ell   \\
\end{array}} \right)
\]
the unknown solution vector and the right-hand side vector of system~(\ref{eq:sysl}), respectively, the solution process with the above multilevel VBARMS factorization consists of a level-by-level forward elimination step followed by an exact solution on the last reduced subsystem and a suitable inverse permutation. The complete solving phase is sketched in Algorithm~\ref{alg:solve}.

\begin{algorithm}
\begin{footnotesize}
 \caption{\texttt{VBARMS\_Solve($A_{\ell+1},b_\ell$)}. The solving phase with the VBARMS method.}\label{alg:solve}
\begin{algorithmic}[1]
\REQUIRE $\ell \in \N^*$, $\ell_{max} \in \N^*$, $b_\ell=\left( f_\ell, g_\ell \right)^T $
\STATE{Solve $L_\ell y = f_\ell$}
\STATE{Compute $g'_\ell=g_\ell-E_\ell U_\ell^{-1}y$}
\IF{$\ell=\ell_{max}$}
\STATE Solve $A_{\ell+1} z_\ell = g_\ell'$ 
\ELSE 
\STATE{Call \texttt{VBARMS\_Solve($A_{\ell+1},g'_\ell$)}}
\ENDIF  
\STATE{Solve $U_\ell y_\ell = \left[ y - L_\ell^{-1} F_\ell z_\ell \right]$}
\end{algorithmic}
\end{footnotesize}
\end{algorithm}

In VBARMS we perform the factorization approximately for memory efficiency. We use block ILU factorization with threshold to invert inexactly both the upper left-most matrix $D_\ell \approx \bar L_\ell \bar U_\ell$, at each level $\ell$, and the last level Schur complement matrix $A_{\ell_{max}} \approx \bar L_S \bar U_S$. The block ILU method used in VBARMS is a straightforward block variant of the one-level pointwise ILUT algorithm. We drop small blocks $B \in \R^{m_B \times n_B}$ in $\bar L_\ell$, $\bar U_\ell$, $\bar L_S$, $\bar U_S$ whenever $\frac{\|B\|_F}{m_B \cdot n_B} < t$, for a given user-defined threshold $t$. The block pivots in block ILU are inverted exactly by using Gaussian Elimination with partial pivoting. Every operation performed during the factorization calls optimized level-3 BLAS routines~\cite{dddh:90}, taking advantage of the finest block structure appearing in the matrices $D_\ell$, $F_\ell$, $E_\ell$, $C_\ell$. Recall that this fine-level block structure results from the block ordering $P_B$ and consists of small, usually dense, blocks in the diagonal blocks of $D_\ell$ as well as in the matrices $E_\ell$, $F_\ell$, $C_\ell$. We do not drop entries in the construction of the Schur complement except at the last level. The same threshold is applied in all these operations.

\begin{algorithm}
\caption{\label{alg:ikj}{\it General ILU Factorization, IKJ Version}.}
\begin{algorithmic}[1]
\REQUIRE A nonzero pattern set $\mathcal{P}$
\FOR{$i=2,\ldots,n$}
  \FOR{$k=1,\ldots,i-1$}
      \IF{$(i,j) \notin \mathcal{P}$} 
		  \STATE $a_{ik}=a_{ik}/a_{kk}$
	  \ENDIF
	  \FOR{$j=k+1,\ldots,n$}	  
	      \IF{$(i,j) \notin \mathcal{P}$} 
			  \STATE $a_{ij}=a_{ij}-a_{ik}a_{kj}$
		  \ENDIF
	  \ENDFOR
  \ENDFOR
\ENDFOR
\end{algorithmic}
\end{algorithm}

\begin{table}
\begin{center}
\begin{tabular}{l|c|c|c|c|c|c|c}
Matrix & Compression & Method & \begin{tabular}{c}Factorization \\ time (s)\end{tabular} & \begin{tabular}{c}Solving \\ time (s)\end{tabular} & \begin{tabular}{c}Total \\ time (s)\end{tabular} & Mem & Its\\ 
\hline
\multirow{6}{*}{HEART1} & \multirow{6}{60px}{\begin{minipage}{60px}\centering Bsize = 18.62 \\
Bdensity = 98.81 \\
$\tau$ = 0.9\end{minipage}} 
    &  &  &  &  &  & \\ 
 &  & VBARMS & 0.12 & 0.43 & 0.55 & 0.83 & 147\\ 
 &  & ILUT & - & - & - & - & -\\ 
 &  & VBILUT & - & - & - & - & -\\ 
 &  & ARMS & - & - & - & - & -\\ 
 &  &  &  &  &  &  & \\ 
\hline
\multirow{6}{*}{PWTK} & \multirow{6}{60px}{\begin{minipage}{60px}\centering Bsize =  56.95 \\
Bdensity = 12.17 \\
$\tau$ = 0.6\end{minipage}} 
    &  &  &  &  &  & \\ 
 &  & VBARMS & 12.71 & 25.02 & 37.73 & 4.42 & 144\\ 
 &  & ILUT & - & - & - & - & -\\ 
 &  & VBILUT & - & - & - & - & -\\ 
 &  & ARMS & - & - & - & - & -\\ 
 &  &  &  &  &  &  & \\ 
\hline
\multirow{6}{*}{RAE} & \multirow{6}{60px}{\begin{minipage}{60px}\centering Bsize = 4.67 \\
Bdensity = 95.83 \\
$\tau$ = 0.8\end{minipage}} 
 &  &  &  &  &  & \\ 
 &  & VBARMS & 1.45 & 1.28 & 2.72 & 2.46 & 34\\ 
 &  & ILUT & - & - & - & - & -\\ 
 &  & VBILUT & - & - & - & - & -\\ 
 &  & ARMS & - & - & - & - & -\\
 &  &  &  &  &  &  & \\   
\hline
\multirow{6}{*}{NASASRB} & \multirow{6}{60px}{\begin{minipage}{60px}\centering Bsize = 9.18 \\
Bdensity = 47.35 \\
$\tau$ = 0.6\end{minipage}} 
 &  &  &  &  &  & \\  
 &  & VBARMS & 2.56 & 3.68 & 6.23 & 3.86 & 76\\ 
 &  & VBILUT & 1.5 & 23.02 & 24.52 & 4.58 & 464\\ 
 &  & ILUT & - & - & - & - & -\\ 
 &  & ARMS & - & - & - & - & -\\ 
 &  &  &  &  &  &  & \\    
\hline
\multirow{6}{*}{OILPAN} & \multirow{6}{60px}{\begin{minipage}{60px}\centering Bsize = 7.01 \\
Bdensity = 99.94 \\
$\tau$ = 0.8\end{minipage}} 
     &        &      &      &      &      & \\ 
 &  & VBARMS & 0.77 & 1.63 & 2.39 & 2.57 & 42\\ 
 &  & ILUT & 0.06 & 32.02 & 32.08 & 0.02 & 952\\ 
 &  & VBILUT & - & - & - & - & -\\ 
 &  & ARMS & - & - & - & - & -\\
 &  &  &  &  &  &  & \\      
\hline
\multirow{6}{*}{BCSSTK35} & \multirow{6}{60px}{\begin{minipage}{60px}\centering Bsize = 11.03 \\
Bdensity = 51.95 \\
$\tau$ = 0.6\end{minipage}} &  &  &  &  &  &  \\ 
 &  & VBILUT & 0.09 & 2.95 & 3.03 & 1.08 & 243\\ 
 &  & VBARMS & 0.15 & 3.22 & 3.36 & 0.95 & 242\\ 
 &  & ILUT & - & - & - & - & -\\ 
 &  & ARMS & - & - & - & - & -\\ 
 &  &  &  &  &  &  & \\      
 \hline 
\end{tabular}
\caption{Assessment performance of VBARMS against other popular preconditioning methods. The symbol `-' means that no convergence is achieved after 1000 iterations of GMRES.}\label{tab:vbarms}
\end{center}
\end{table}

\ms

\subsection{The new implementation of VBARMS}
The code for the VBARMS method is developed in the C language and is adapted from the existing ARMS code available in the ITSOL package~\cite{itsol:09}. The compressed sparse storage format of ARMS is modified to store block vectors and block matrices of variable size as a collection of contiguous nonzero dense blocks (we refer to this data storage format as VBCSR). However, the implementation used in this paper is different and noticeably faster than the one described in~\cite{VBARMS}.
In the old implementation, the approximate transformation matrices $E_\ell \bar U_\ell^{-1}$ and $\bar L_\ell^{-1} F_l$ appearing in Eqn~(\ref{eq:sysl}) at step $\ell$ were explicitly computed and temporarily stored in the VBCSR format. They were discarded from the memory immediately after assembling $A_{\ell+1}$. In the new implementation, we first compute the factors $\bar L_\ell$, $\bar U_\ell$ and $\bar L_\ell^{-1} F_\ell$ by performing a variant of the IKJ version of the Gaussian Elimination algorithm (Algorithm~\ref{alg:ikj}), where index $I$ runs from $2$ to $m_\ell$, index $K$ from $1$ to $(I-1)$ and index $J$ from $(K+1)$ to $n_\ell$.
This loop applies implicitly $\bar L_\ell^{-1}$ to the block row 
$\left[ {D_\ell ~,~ F_\ell  } \right]$ to produce $\left[ {U_\ell ~,~ \bar L_\ell ^{ - 1} F_\ell  } \right]$. In the second loop, Gaussian Elimination is performed on the block row $\left[ {E_\ell ~,~ C_\ell  } \right]$ using the multipliers computed in the first loop to give ${E_\ell  \bar U_\ell ^{ - 1}}$ and an approximation of the Schur complement $A_{\ell+1}$. We explicitly permute the matrix after Step~1 at the first level as well as the matrices involved in the factorization at each new reordering step. The improvement of efficiency obtained with the new implementation is noticeable, as appears from the results shown in Table~\ref{tab:oldnew}. Finally, in Table~\ref{tab:vbarms} we assess the performance of the VBARMS method against other popular preconditioning techniques;
we report on the number of GMRES iterations required to reduce the initial residual by $6$ orders of magnitude using a block incomplete LU factorization as a preconditioner for GMRES.
The results show a remarkable robustness for low to moderate memory cost.

\begin{table}
\begin{center}
\begin{tabular}{p{50px}|c|c|c|c|c|c}
Matrix & Implementation & \begin{tabular}{c}Factorization \\ time (s)\end{tabular} & \begin{tabular}{c}Solving \\ time (s)\end{tabular} & \begin{tabular}{c}Total \\ time (s)\end{tabular} & Mem & Its\\ 
\hline
\multirow{2}{*}{HEART1}  & New & 0.12 & 0.43 & 0.55 & 0.83 & 147\\ 
 & Old & 0.36 & 0.33 & 0.69 & 0.86 & 113\\ 
\hline
\multirow{2}{*}{PWTK} & New & 12.71 & 25.02 & 37.73 & 4.42 & 144\\ 
 & Old & 90.73 & 26.08 & 116.81 & 4.95 & 140\\ 
\hline
\multirow{2}{*}{RAE} & New & 1.45 & 1.28 & 2.72 & 2.46 & 34\\ 
 & Old & 5.12 & 1.15 & 6.27 & 2.71 & 30\\ 
\hline
\multirow{2}{*}{NASASRB} & New & 2.56 & 3.68 & 6.23 & 3.86 & 76\\ 
 & Old & 15.54 & 3.34 & 18.88 & 4.06 & 64\\ 
\hline
\multirow{2}{*}{OILPAN} & New & 0.77 & 1.63 & 2.39 & 2.57 & 42\\ 
 & Old & 5.64 & 1.29 & 6.93 & 2.62 & 32\\ 
\hline
\multirow{2}{*}{BCSSTK35} & New & 0.15 & 3.22 & 3.36 & 0.95 & 242\\ 
 & Old & - & - & - & - & -\\ \hline
\end{tabular}
\caption{Comparative experiments with the old and the new VBARMS codes, implementing a different partial (block) factorization step.  The symbol `-' means that no convergence is achieved after 1000 iterations of GMRES.}\label{tab:oldnew}
\end{center}
\end{table}

\section{Using VBARMS in parallel computing}~\label{sec:4}

\begin{figure}[ht!]
\begin{center}\includegraphics[width=7.0cm]{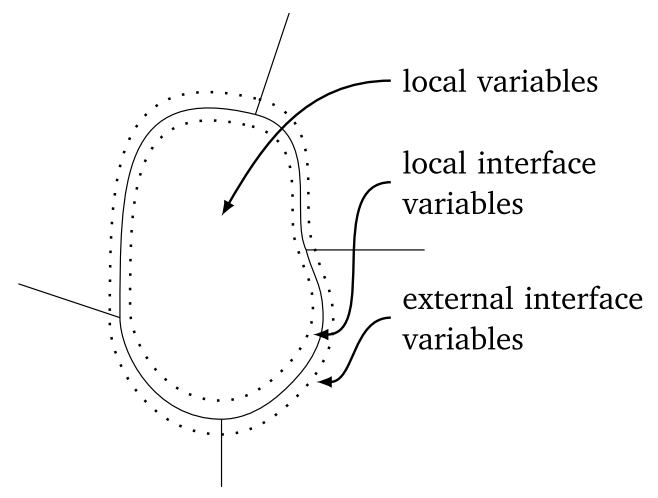}
\caption{Local domain from a physical viewpoint.}
\label{fig:distributed}
\end{center}
\end{figure}

In the experiments reported in this section the VBARMS method is used for solving large linear systems on distributed memory computers; its overall performance are assessed against the parallel implementation of the ARMS solver provided in the pARMS package~\cite{pARMS}. On multicore machines, the quotient graph $\mathcal{G}/\mathcal{B}$ is split into distinct subdomains using a parallel graph partitioner, and each of them is assigned to a different core. We follow the parallel framework described in~\cite{pARMS} which separates the nodes assigned to the $i$th subdomain into {\it interior nodes}, that are those coupled only with local variables by the equations, and {\it interface nodes}, those that may be coupled with local variables stored on processor $i$ as well as with remote variables stored on other processors (see~Figure~\ref{fig:distributed}). 
The vector of the local unknowns $x_i$ and the local right-hand side $b_i$ are split accordingly in two separate components: the subvector corresponding to the internal nodes followed by the subvector of the local interface variables
\[
x_i  = \left( {\begin{array}{*{20}c}
   {u_i }  \\
   {y_i }  \\
\end{array}} \right),~~~ b_i  = \left( {\begin{array}{*{20}c}
   {f_i }  \\
   {g_i }  \\
\end{array}} \right).
\]
The rows of $A$ corresponding to the nodes belonging to the $i$th subdomain are assigned to the $i$th processor. They are naturally separated into a local matrix $A_i$ acting on the local variables $x_i=(u_i,y_i)^T$, and an interface matrix $U_i$ acting on the remotely stored subvectors of the external interface variables $y_{i,ext}$. Hence we can write the local equations on processor $i$ as
\[
 A_i x_i  + U_{i,ext} y_{i,ext} = b_i
\]
or, in expanded form, as 
\begin{equation}\label{eq:paral}
\left( {\begin{array}{*{20}c}
   {B_i } & {F_i }  \\
   {E_i } & {C_i }  \\
\end{array}} \right)\left( {\begin{array}{*{20}c}
   {u_i }  \\
   {y_i }  \\
\end{array}} \right) + \left( {\begin{array}{*{20}c}
   0  \\
   {\sum\nolimits_{j \in N_i } {E_{ij} y_j } }  \\
\end{array}} \right) = \left( {\begin{array}{*{20}c}
   {f_i }  \\
   {g_i }  \\
\end{array}} \right),
\end{equation}
where $N_i$ is the set of subdomains that are neighbors to subdomain $i$ and the submatrix $E_{ij} y_j$ accounts for the contribution to the local equation from the $j$th neighboring subdomain. Notice that matrices $B_i$, $C_i$, $E_i$ and $F_i$ still preserve the finest block structure imposed by the block ordering $P_B$. At this stage, the VBARMS method described in Section~\ref{sec:3} can be used as a local solver for different types of global preconditioners.

\ms 

In the simplest parallel implementation, the so-called block-Jacobi preconditioner, the sequential VBARMS method can be applied to invert approximately each local matrix $A_i$. The standard Jacobi iteration for solving $Ax=b$ is defined as
\[
x_{n + 1}  = x_n  + D^{ - 1} \left( {b - Ax_n } \right) = D^{ - 1} \left( {Nx_n  + b} \right)
\]
where $D$ is the diagonal of $A$, $N=D-A$ and $x_0$ is some initial approximation. In cases we have a graph partitioned matrix, the matrix $D$ is block diagonal and the diagonal blocks of $D$ are the local matrices $A_i$. The interest to consider this basic approach is its inherent parallelism, since the solves with the matrices $A_i$ are performed independently on all the processors and no communication is required.

\ms 

If the diagonal blocks of the matrix $D$ are enlarged in the block-Jacobi method so that they overlap slightly, the resulting preconditioner is called Schwarz preconditioner. Consider again a graph partitioned matrix with $N$ nonoverlapping sets $W_i^0$, $i=1,\ldots,N$ and $W_0  = \bigcup\nolimits_{i = 1}^N {W_0^i}$. We define a $\delta$-overlap partition
\[
W^\delta   = \bigcup\limits_{i = 1}^N {W_i^\delta}
\]
where $W_i^\delta   = adj\left( {W_i^{\delta  - 1} } \right)$
and $\delta > 0$ is the level of overlap with the neighbouring domains. For each subdomain, we define a restriction operator $R^\delta_i$, which is an $n \times n$ matrix with the $(j,j)$th element equal to $1$ if $j\in W_i^\delta$, and zero elsewhere. We then denote
\[
A_i  = R_i^\delta  AR_i^\delta .
\]
The global preconditioning matrix $M_{RAS}$ is defined as
\[
M^{ - 1}_{RAS}  = \sum\limits_{i = 1}^s {R_i^T A_i^{ - 1} R_i } .
\]
and named as the Restricted Additive Schwarz preconditioner (RAS)~\cite{quva:99,SAAD-SHORT}. Note that the preconditioning step is still parallel, as the different components of the error update are formed independently. However, some communication is required in the final update, as the components are added up from each subdomain due to overlapping.

\bs 

A third global preconditioner that we consider in this study is based on the Schur complement approach. 
In Eqn~(\ref{eq:paral}), we can eliminate the vector of interior unknowns $u_i$ from the first equations to compute the local Schur complement system
\[
S_i y_i  + \sum\limits_{j \in N_i } {E_{ij} } y_j  = g_i  - E_i B_i^{ - 1} f_i  \equiv g'_i ,
\]
where $S_i$ denotes the local Schur complement matrix
\[
	S_i = C_i - E_i B_i^{-1} F_i.
\]
The local Schur complement equations considered altogether write as the global Schur complement system
\begin{equation}\label{eq:gschur}
\left( {\begin{array}{*{20}c}
   {S_1 } & {E_{12} } &  \ldots  & {E_{1p} }  \\
   {E_{21} } & {S_2 } &  \ldots  & {E_{2p} }  \\
    \vdots  & {} &  \ddots  &  \vdots   \\
   {E_{p1} } & {E_{p - 1,2} } &  \ldots  & {S_p }  \\
\end{array}} \right)\left( {\begin{array}{*{20}c}
   {y_1 }  \\
   {y_2 }  \\
    \vdots   \\
   {y_p }  \\
\end{array}} \right) = \left( {\begin{array}{*{20}c}
   {g'_1 }  \\
   {g'_2 }  \\
    \vdots   \\
   {g'_p }  \\
\end{array}} \right) ,
\end{equation}
where the off-diagonal matrices $E_{ij}$ are available from the parallel distribution of the linear system. One preconditioning step with the Schur complement preconditioner consists in solving approximately the global system~(\ref{eq:gschur}), and then recovering the $u_i$ variables from the local equations as
\begin{equation}\label{loc_var_schur}
	u_i=B_i^{-1}[f_i-F_i y_i]
\end{equation}
at the cost of one local solve. We solve the global system~(\ref{eq:gschur}) by running a few steps of the GMRES method preconditioned by a block diagonal matrix, where the diagonal blocks are the local Schur complements $S_i$. The factorization 
\[
	S_i = L_{S_i } U_{S_i }
\]	
is obtained as by-product of the LU factorization of the local matrix $A_i$,
\[
A_i  = \left( {\begin{array}{*{20}c}
   {L_{B_i } } & 0  \\
   {E_i U_{B_i }^{ - 1} } & {L_{S_i } }  \\
\end{array}} \right)\left( {\begin{array}{*{20}c}
   {U_{B_i } } & {L_{B_i }^{ - 1} F_i }  \\
   0 & {U_{S_i } }  \\
\end{array}}  \right).
\]
which is by the way required to compute the $u_i$ variables in~(\ref{loc_var_schur}).

\subsection{Experiments}

Some preliminary results with a parallel MPI-based implementation of VBARMS for distributed memory computers, reported in a conference contribution~\cite{ppam2014}, revealed promising performance against the parallel ARMS method and the conventional ILUT method. They showed that exposing dense matrix blocks during the factorization may lead to more efficient and more stable parallel solvers. The parallel implementation of VBARMS considered in this study differs from the one presented in~\cite{ppam2014} in one important aspect. In the old implementation we used a sequential graph partitioner, namely the recursive dissection partitioner from the METIS package~\cite{metis4}, to split the quotient graph $G/{\mathcal{B}}$ and then assign the computed partitions to different processors. 
In the new implementation, the quotient graph is initially distributed amongst the available processors; then, the built-in parallel hypergraph partitioner available in the Zoltan package~\cite{ZoltanHomePage} is applied on the distributed data structure to compute an optimal partitioning of the quotient graph that can minimize the amount of communications.

In the experiments reported in Table~\ref{tab:zoltan_sven_16} we notice the significant reduction of CPU time spent for the graph partitioning operation in the new implementation of VBARMS; note that the numerical efficiency of the solvers is generally well preserved. The matrix problems used are listed in Table~ \ref{tab:matrices_par}. The parallel experiments were run on the large-memory nodes (32 cores/node and 1TB of memory) of the TACC Stampede system located at the University of Texas at Austin. TACC Stampede is a 10 PFLOPS (PF) Dell Linux Cluster based on 6,400+ Dell PowerEdge server nodes, each outfitted with 2 Intel Xeon E5 (Sandy Bridge) processors and an Intel Xeon Phi Coprocessor (MIC Architecture). We linked the vendor BLAS library on Stampede, which has BLAS via MKL loaded by default and is multi-threaded. We used the Flexible GMRES (FGMRES) method~\cite{saad:93} as Krylov subspace method, a tolerance of $1.0e-6$ in the stopping criterion and a maximum number of iteration equal to $1000$. Memory costs were calculated as the ratio between the sum of the number of nonzeros in the local preconditioners, and the sum of the number of nonzeros in the local matrices $A_i$. Overall, the Restricted Additive Schwarz solver showed better performance against the Block Jacobi and the Schur-complement methods.

\begin{table}[!ht]
 \begin{center}
  \begin{tabular}{l|c|l|r}
    \hline
    Name & Size & Application & nnz(A) \\
    \hline
	AUDIKW\_1 & 943695 & Structural problem & 77651847 \\
	LDOOR & 952203 & Structural problem & 42493817 \\
	STA004 & 891815 & Fluid Dynamics & 55902989 \\
	STA004 & 891815 & Fluid Dynamics & 55902989 \\
	\hline
  \end{tabular}
 \caption{\label{tab:matrices_par}Set and characteristics of test matrix problems.}
 \end{center}
\end{table}

\begin{table}[!h]
\begin{center}
\begin{footnotesize}
\begin{tabular}{c|c|c|c|c|c|c|c|c}
Matrix  & Method  & \begin{tabular}{c}Graph \\ type (s)\end{tabular} & \begin{tabular}{c}Graph \\ time (s)\end{tabular} & \begin{tabular}{c}Factorization \\ time (s)\end{tabular} & \begin{tabular}{c}Solving \\ time (s)\end{tabular} & \begin{tabular}{c}Total \\ time (s)\end{tabular} & Its &  Mem \\
\hline
{AUDIKW$\_$1} &
\tabincell{c}{BJ+VBARMS\\\\RAS+VBARMS\\\\SCHUR+VBARMS}&
\tabincell{c}{METIS (seq.)\\Zoltan (par.)\\METIS (seq.)\\Zoltan (par.)\\METIS (seq.)\\Zoltan (par.)}&
\tabincell{c}{54.5 \\ 5.2\\54.2\\5.3 \\ 54.4 \\ 5.3}&
\tabincell{c}{18.88\\ 17.28\\ 19.54\\ 22.75 \\ 82.72 \\ 166.09} &
\tabincell{c}{51.35\\ 37.98\\ 26.68\\ 22.24 \\ 295.11 \\ 327.06} &
\tabincell{c}{70.23\\ 55.26\\ 46.22\\ 44.99 \\ 377.83 \\ 493.15} &
\tabincell{c}{136\\117\\ 46\\52 \\ 69 \\ 59} &
\tabincell{c}{3.13\\ 2.74\\2.93 \\2.87 \\ 6.21 \\ 4.60}\\
\hline
{LDOOR} &
\tabincell{c}{BJ+VBARMS\\\\RAS+VBARMS\\\\SCHUR+VBARMS}&
\tabincell{c}{METIS (seq.)\\Zoltan (par.)\\METIS (seq.)\\Zoltan (par.)\\METIS (seq.)\\Zoltan (par.)}&
\tabincell{c}{30.0 \\ 1.1 \\ 29.0 \\ 1.1 \\ 29.0 \\ 1.1}&
\tabincell{c}{1.29\\ 1.04\\1.56\\1.12 \\ 5.81 \\ 5.64} &
\tabincell{c}{25.10\\ 18.09\\13.40\\ 12.73 \\ 16.75 \\ 4.78} &
\tabincell{c}{26.40\\ 19.12\\14.95\\ 13.85 \\ 22.56 \\ 10.42} &
\tabincell{c}{345\\273 \\ 200\\196 \\ 54 \\ 37} &
\tabincell{c}{1.95\\ 1.95\\ 2.00\\1.99 \\ 3.63 \\ 3.32}\\
\hline
{STA004} &
\tabincell{c}{BJ+VBARMS\\\\RAS+VBARMS\\\\SCHUR+VBARMS}&
\tabincell{c}{METIS (seq.)\\Zoltan (par.)\\METIS (seq.)\\Zoltan (par.)\\METIS (seq.)\\Zoltan (par.)}&
\tabincell{c}{79.4 \\ 2.5 \\ 81.7 \\ 2.6 \\ 81.4 \\ 2.5}&
\tabincell{c}{7.53\\ 5.11\\ 9.55\\7.90 \\ 17.46 \\ 16.05} &
\tabincell{c}{42.56\\ 24.12\\ 34.27\\ 23.09 \\ 135.58 \\ 113.24} &
\tabincell{c}{50.08\\ 29.23\\ 43.82\\30.99 \\153.04 \\ 129.28} &
\tabincell{c}{90\\ 72 \\42\\34 \\90 \\ 88} &
\tabincell{c}{3.61\\ 3.61\\3.85\\3.31 \\ 5.29 \\ 5.40}\\
\hline
{STA008} &
\tabincell{c}{BJ+VBARMS\\\\RAS+VBARMS\\\\SCHUR+VBARMS}&
\tabincell{c}{METIS (seq.)\\Zoltan (par.)\\METIS (seq.)\\Zoltan (par.)\\METIS (seq.)\\Zoltan (par.)}&
\tabincell{c}{81.9 \\ 2.3 \\ 81.8 \\ 2.4 \\ 81.2 \\ 2.4 }&
\tabincell{c}{11.36\\ 9.45\\ 15.01 \\12.90 \\ 56.20 \\ 66.42} &
\tabincell{c}{85.77\\ 50.17\\ 67.98\\ 46.52 \\ 564.75 \\ 490.25} &
\tabincell{c}{97.14\\ 59.62\\ 82.99\\ 59.42 \\ 620.94 \\556.67} &
\tabincell{c}{227\\ 170\\ 101\\97 \\ 188 \\ 201} &
\tabincell{c}{4.77\\ 4.78\\ 5.10\\ 5.07 \\ 8.94 \\ 9.83}\\
\hline
\end{tabular}
\end{footnotesize}
\end{center}
\caption{Performance comparison of serial and parallel graph partition on 16 processors.
Notation: P-N means number of processors, G-Type means graph partitioning strategy, G-time means 
partitioning timing cost, P-T means preconditioning construction time, I-T iterative solution time, Mem
means memory costs.}\label{tab:zoltan_sven_16}
\end{table}

\FloatBarrier

\subsection{A case study in large-scale turbulent flows analysis}

We finally get back to the starting point that motivated this study. In this section we present a performance analysis with the parallel VBARMS implementation for solving large block structured linear systems arising from an implicit Newton-Krylov formulation of the Reynolds Averaged Navier Stokes (briefly, RANS) equations. 
Although explicit multigrid techniques have dominated the Computational Fluid Dynamics (CFD) arena for a long time, implicit methods based on Newton's rootfinding algorithm are recently receiving increasing attention because of their potential to converge in a very small number of iterations. 
One of the most recent outstanding examples on the use of implicit unstructured RANS CFD is provided in the article~\cite{wozi:08}, which reports the turbulent analysis of the flow past three-dimensional wings using a vertex-based unstructured Newton-Krylov solvers. Practical implicit CFD solvers need to be combined with ad-hoc preconditioners to invert efficiently the large nonsymmetric linear system at each step of Newton's algorithm. 

\ms

\begin{figure}[ht!]
\begin{subfigmatrix}{2}
\subfigure[The flux balance of cell $T$ is scattered among its vertices.]%
  {\includegraphics[width=6.0cm]{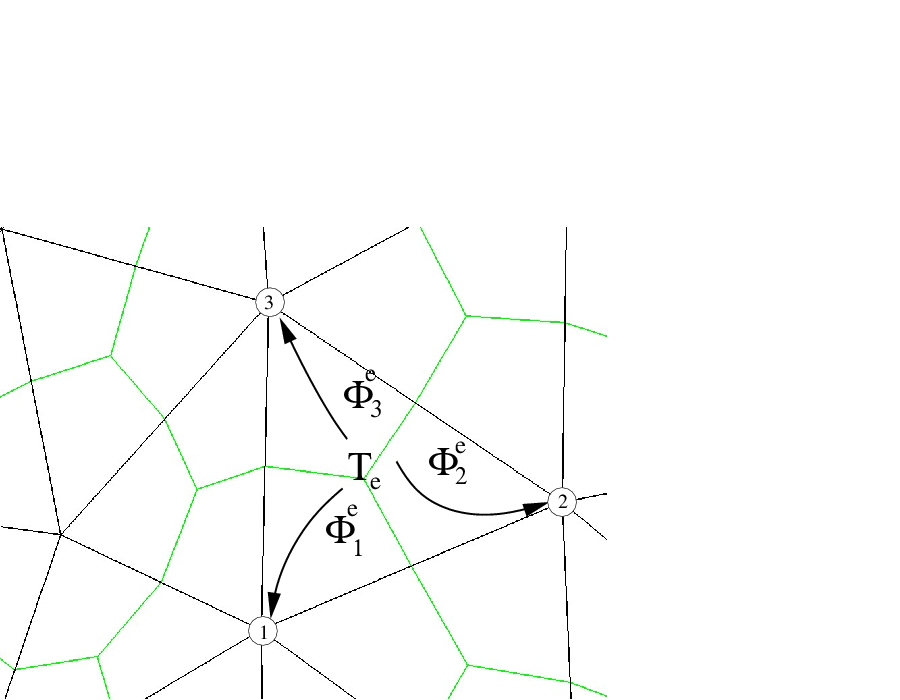}\label{fig:scatter}}
\subfigure[Gridpoint $i$ gathers the fractions of cell residuals from the surrounding cells.]%
  {\includegraphics[width=6.0cm]{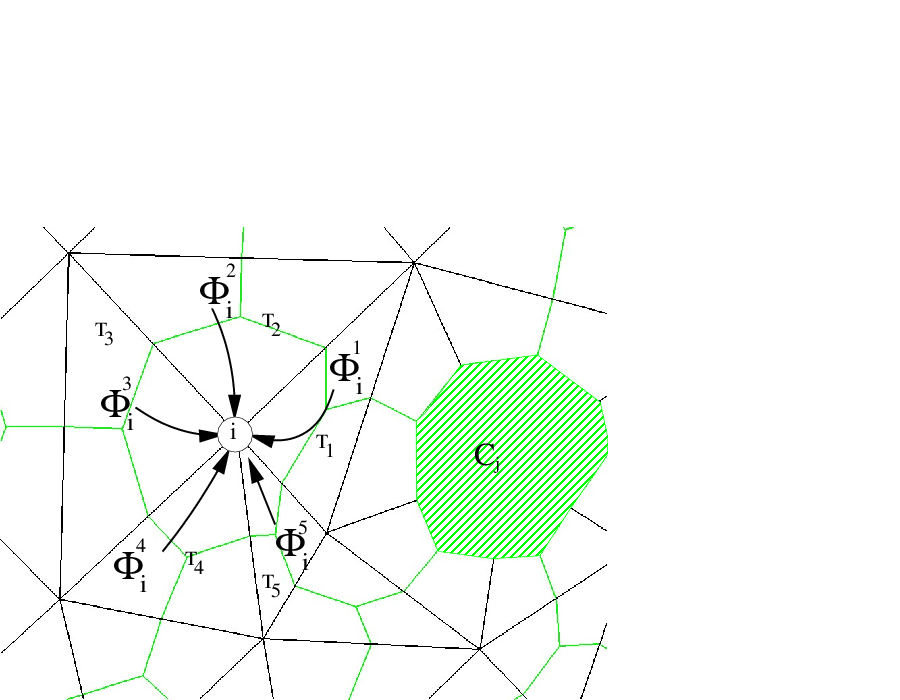}\label{fig:distr}}
\end{subfigmatrix}
  \caption{Residual distribution concept.}
\label{fig:concept}
\end{figure}

Throughout this section we use standard notation for the kinematic and thermodynamic variables: we denote by $\vec{u}$ the flow velocity, by $\rho$ the density, $p$ is the pressure, $T$ is the temperature, $e$ and $h$ are respectively the specific total energy and enthalpy, ${\nu}$ is the laminar kinematic viscosity and $\tilde{\nu}$ is a scalar variable related to the turbulent eddy viscosity via a damping function. The quantity $a$ denotes the sound speed or the square root of the artificial compressibility constant in case of the compressible, respectively incompressible, flow equations. 
In the case of high Reynolds number flows, we account for turbulence effects by the 
RANS equations that are obtained from the Navier-Stokes (NS) equations by means of a
time averaging procedure. The RANS equations have the same structure as the NS equations 
with an additional term, the Reynolds' stress tensor, that accounts for the effects of 
the turbulent scales on the mean field. Using Boussinesq's approximation, the 
Reynolds' stress tensor is linked to the mean velocity gradient through the turbulent (or eddy) viscosity. In our study, the turbulent viscosity is modeled using the Spalart-Allmaras one-equation model~\cite{spal:94}. The physical domain is partitioned into nonoverlapping control volumes drawn around each gridpoint by joining, in two space dimensions, the centroids of gravity of the surrounding cells with the midpoints of all the edges that connect that gridpoint with its nearest neighbors, as shown in Figure~\ref{fig:concept}. 

Given a control volume $C_i$, fixed in space and bounded by the control
surface $\partial C_i$ with inward normal ${\vec n}$, we write the
governing conservation laws of mass, momentum, energy and turbulence
transport equations as
\begin{equation} \label{Carpentieri_contrib_discrete_compressible}
   \int_{C_i} {\frac{\partial {{\vec q}_i}}{\partial {t}}} \, dV =
     \oint_{\partial C_i} {\vec n} \cdot {\vec F} \, dS
   - \oint_{\partial C_i} {\vec n} \cdot {\vec G} \, dS +
   \int_{C_i} {\vec s} \, dV,
\end{equation} where we denote by ${\vec q}$ the vector of conserved variables. For
compressible flows, we have ${\vec q}=\left( \rho, \rho e , \rho {\vec u} , 
\tilde{\nu} \right)^T,$ and for incompressible, constant
density flows,  ${\vec q}=\left( p ,  {\vec u} , \tilde{\nu}
\right)^T.$ In (\ref{Carpentieri_contrib_discrete_compressible}), the vector operators
${\vec F}$ and ${\vec G}$ represent the inviscid and viscous
fluxes, respectively. For compressible flows, we have
\[
\label{Carpentieri_contrib_compr_invscflx}
   {\vec F} = \left( \begin{array}{c}
   \rho \vec{u} \\ \rho \vec{u} h \\ \rho \vec{u}   \vec{u} + p {\bf I} \\
   \tilde{\nu} \vec{u}
   \end{array} \right), \quad
   {\vec G} = \frac{1}{{\operatorname{Re} _\infty  }} \left( \begin{array}{c}
   0 \\ \vec{u} \cdot { {{{\tau}}} } + \nabla q \\ {{{\tau}}} \\
\frac{1}{\sigma} \left[ \left( {{\nu}} + \tilde{\nu} \right) \nabla
\tilde{\nu} \right]
   \end{array} \right),
\]
\noindent and for incompressible, constant density flows,
\[
{\vec F} = \left( \begin{array}{c}
  a^2 \vec{u}  \\
  \vec{u}\vec{u} + p {\bf I}  \\
  \tilde{\nu} \vec{u} \\
\end{array}  \right),~~{\vec G} = \frac{1}
{{\operatorname{Re} _\infty  }}\left( \begin{array}{c}
  0  \\
  {\tau}  \\
\frac{1}{\sigma} \left[ \left( {\nu} + \tilde{\nu} \right) \nabla
\tilde{\nu} \right]
\end{array}  \right),
\]
where $\tau$ is the Newtonian stress tensor.
The source term vector $\vec{s}$ has a non-zero entry only in the row
corresponding to the turbulence transport equation, which takes the form
\begin{equation}
\label{Carpentieri_contrib_source_term}
 c_{b1} \left[ 1 - f_{t2} \right] \tilde{S} \tilde{\nu}
+ \frac{1}{\sigma Re} \left[
c_{b2} \left( \nabla \tilde{\nu} \right)^2 \right] + \\
- \frac{1}{Re} \left[ c_{w1} f_{w} - \frac{c_{b1}}{\kappa^2} f_{t2}
\right] \left[ \frac{\tilde{\nu}}{d} \right]^2.   
\end{equation}

For a description of the various functions and constants involved 
in~(\ref{Carpentieri_contrib_source_term}) we refer the reader
to~\cite{spal:94}. 

\ms

We consider a fluctuation splitting approach to discretize in space 
the integral form of the governing equations~(\ref{Carpentieri_contrib_discrete_compressible})
over each control volume $C_i$. The flux integral is evaluated over 
each triangle (or tetrahedron) in the mesh, and then split among 
its vertices~\cite{ijcfd:aldo} (see Figure~\ref{fig:concept}), so that 
we may write from Eq.~(\ref{Carpentieri_contrib_discrete_compressible})
$$
\label{Carpentieri_contrib:_conservationlaw_ib}
   \int_{C_i} {{\frac{\partial {{\vec q}_i}}{\partial {t}}}}  \, dV =
   \sum_{T \ni i} {\vec \phi}_i^T
$$ where
\[
  {\vec \phi}^T = \oint_{\partial T} {\vec n} \cdot {\vec F} \, dS
   - \oint_{\partial T} {\vec n} \cdot {\vec G} \, dS +
   \int_{T} {\vec s} \, dV
\]
is the flux balance evaluated over cell $T$ and $\vec{\phi}_i^T$ is the
fraction of cell residual scattered to vertex $i$. 
Upon discretization of the governing equations, we obtain a system
of ordinary differential equations of the form
\begin{equation} \label{Carpentieri_contrib_semi_d_matrix}
    {M} \frac{d {\vec q}}{dt} = {\vec r}({\vec q}),
\end{equation} where $t$ denotes the pseudo time variable, ${M}$ is the mass
matrix and ${\vec r}({\vec q})$ represents the nodal residual vector of spatial
discretization operator, which vanishes at steady state. The
residual vector is a (block) array of dimension equal to the number
of meshpoints times the number of dependent variables, $m$; for a
one-equation turbulence model, $m=d+3$ for compressible flows and
$m=d+2$ for incompressible flows, $d$ being the spatial dimension.
If the time derivative in equation (\ref{Carpentieri_contrib_semi_d_matrix}) is
approximated using a two-point one-sided finite difference (FD) formula
we obtain the following implicit scheme:
\begin{equation} \label{Carpentieri_contrib_big_system}
   \left( \frac{1}{\Delta t^n}{V} - {J} \right) \left(
   {\vec q}^{n+1} - {\vec q}^n \right) = {\vec r}({\vec q}^n),
\end{equation}
where we denote by ${J}$ the Jacobian of the residual
$\displaystyle {{\frac{\partial {{\vec r}}}{\partial {{\vec q}}}}} $. 
We use a finite difference approximation of the Jacobian, where the individual entries
of the vector of nodal unknowns are perturbed by a small amount
$\epsilon$ and the nodal residual is then recomputed for the
perturbed state. Eq. (\ref{Carpentieri_contrib_big_system}) represents a
large nonsymmetric sparse linear
system of equations to be solved at each pseudo-time step for the
update of the vector of the conserved variables. 
The nonzero pattern of the sparse coefficient matrix is 
symmetric; on average, the number of non-zero (block) entries per row 
in our discretization scheme equals 7 in 2D and 14 in 3D. 
Choice of the iterative solver and of the preconditioner can have a strong influence 
on computational efficiency, especially when the mean flow and turbulence transport 
equations are solved in fully coupled form like we do. 

\ms 

\begin{figure}[ht!]
\begin{center}
\begin{minipage}{4cm}
\begin{center}
  {\includegraphics[width=4.0cm]{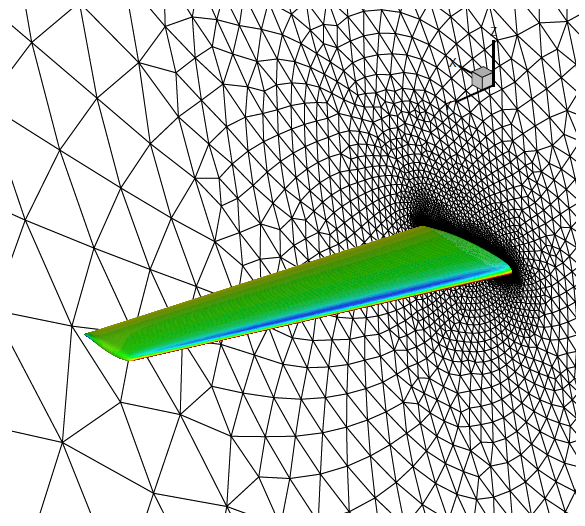}}
\end{center}
\end{minipage}
\hfil
\begin{minipage}{6cm}
\begin{center}
\begin{tabular}{lll}
Ref. Area, &S = 290322 mm$^2$ &= 450 in$^2$\\
Ref. Chord,&c = 197.556 mm &= 7.778 in\\
Ref. Span, &b = 1524 mm &= 60 in
\end{tabular}

\medskip 

\begin{tabular}{lrr}
RANS1 : & $n = ~4918165$ & $nnz = ~318370485$ \\  
RANS2 : & $n = ~4918165$ & $nnz = ~318370485$ \\  
RANS3 : & $n = ~9032110$ & $nnz = ~670075950$ \\  
RANS4 : & $n = 12085410$ & $nnz = ~893964000$ \\  
RANS5 : & $n = 22384845$ & $nnz = 1659721325$  

\end{tabular}
\end{center}
\end{minipage}
\end{center}
\caption{Geometry and mesh characteristics of the DPW3 Wing-1 problem proposed in the 3rd AIAA Drag Prediction Workshop. Note that problems RANS1 and RANS2 correspond to the same mesh, and are generated at two different Newton steps.}\label{fig:dpw3}
\end{figure}

We consider turbulent incompressible flow analysis past a three-dimensional wing
illustrated in Fig.~\ref{fig:dpw3}. 
The geometry, called DPW3 Wing-1, was proposed in the 3rd AIAA Drag Prediction Workshop~\cite{drag}. 
Flow conditions are 0.5$^{\circ}$ angle of attack and Reynolds number based on the reference chord equal to $5 \cdot 10^6$. The freestream turbulent viscosity is set to 10\% of its laminar value.

In Table~\ref{tab:cfd_results} we show experiments with parallel VBARMS on the first four meshes of the DPW3 Wing-1 problem. We illustrate only examples with the parallel graph partitioning strategy described in Section~\ref{sec:4}. In Table~\ref{tab:cfd_results_RANS5} we report on only one experiment on the largest mesh, as this is a resource demanding problem. In Table~\ref{tab:cfd_zoltan} we perform a strong scalability study on the problem denoted as RANS2 by increasing the number of processors. 
Finally, in Table~\ref{tab:cfd_arms} we report on comparative results with parallel VBARMS against other popular solvers. The method denoted as \texttt{pARMS} is the solver described in~\cite{pARMS}, using default parameters. The results of our experiments confirm the same trend of performance shown on general problems. The proposed VBARMS method is remarkably efficient for solving block structured linear systems arising in applications in combination with conventional parallel global solvers such as in particular the Restricted Additive Schwarz preconditioner. A truly parallel implementation of the VBARMS method that may offer better numerical scalability will be considered as the next step of this research.

\begin{table}[!h]
\begin{center}
\begin{footnotesize}
\begin{tabular}{c|c|c|c|c|c|c|c}
Matrix  & Method & \begin{tabular}{c}Graph \\ time (s)\end{tabular} & \begin{tabular}{c}Factorization \\ time (s)\end{tabular} & \begin{tabular}{c}Solving \\ time (s)\end{tabular} & \begin{tabular}{c}Total \\ time (s)\end{tabular} & Its &  Mem \\
\hline

RANS1 &
\tabincell{c}{BJ+VBARMS\\RAS+VBARMS\\SCHUR+VBARMS}&
\tabincell{c}{17.3\\ 17.4 \\17.6}&
\tabincell{c}{8.58\\ 10.08 \\ 11.94} &
\tabincell{c}{41.54\\ 42.28 \\ 55.99} &
\tabincell{c}{50.13\\ 52.37 \\ 67.93} &
\tabincell{c}{34\\ 19 \\35} &
\tabincell{c}{2.98\\ 3.06 \\2.57}\\
\hline
RANS2 &
\tabincell{c}{BJ+VBARMS\\RAS+VBARMS\\SCHUR+VBARMS}&
\tabincell{c}{17.0 \\ 16.8 \\ 17.5}&
\tabincell{c}{16.72\\ 21.65 \\ 168.85} &
\tabincell{c}{70.14\\ 80.24 \\ 173.54} &
\tabincell{c}{86.86\\ 101.89 \\ 342.39} &
\tabincell{c}{47 \\ 39 \\ 24} &
\tabincell{c}{4.35\\ 4.49 \\ 6.47}\\
\hline
RANS3 &
\tabincell{c}{BJ+VBARMS\\RAS+VBARMS\\SCHUR+VBARMS}&
\tabincell{c}{27.2 \\ 25.2 \\ 22.0}&
\tabincell{c}{99.41\\ 119.32 \\ 52.65} &
\tabincell{c}{187.95\\ 90.47\\ 721.67} &
\tabincell{c}{287.36\\ 209.79 \\ 774.31} &
\tabincell{c}{154\\ 71 \\140} &
\tabincell{c}{4.40\\ 4.48 \\4.39}\\
\hline
RANS4 &
\tabincell{c}{BJ+VBARMS\\RAS+VBARMS\\SCHUR+VBARMS}&
\tabincell{c}{51.5 \\ 43.9 \\ 39.3}&
\tabincell{c}{12.05\\ 14.05\\ 15.14} &
\tabincell{c}{105.89\\ 91.53 \\ 289.89} &
\tabincell{c}{117.94\\ 105.58\\ 305.03} &
\tabincell{c}{223 \\ 143\\ 179} &
\tabincell{c}{3.91\\ 4.12\\3.76}\\
\hline
\end{tabular}
\end{footnotesize}
\end{center}
\caption{Experiments on the DPW3 Wing-1 problem. The RANS1, RANS2 and RANS3 test cases are solved on 32 processors, whereas the RANS4 problem on 128 processors.}\label{tab:cfd_results}
\end{table}

\begin{table}[!ht]
\begin{center}
\begin{tabular}{c|c|c|c|c}
Matrix  & Method & Total time (s) & Its &  Mem \\
\hline
\tabincell{c}{RANS5}&
\tabincell{c}{RAS+VBARMS\\}&
\tabincell{c}{291.42\\} &
\tabincell{c}{235} &
\tabincell{c}{4.05}\\
\hline

\end{tabular}
\end{center}
\caption{RANS5 problem is solved on 128 processors.}\label{tab:cfd_results_RANS5}
\end{table}

\begin{table}
\begin{center}
\begin{tabular}{c|c|c|c|c|c}
Solver  & \begin{tabular}{c}Number of \\ processors \end{tabular}  & \begin{tabular}{c}Graph \\ time (s)\end{tabular}  & \begin{tabular}{c}Total \\ time (s)\end{tabular} & Its &  Mem \\
\hline
RAS+VBARMS & 
\tabincell{c}{ 8\\ 16\\ 32\\ 64\\ 128}& 
\tabincell{c}{ 38.9\\ 28.0\\ 17.0\\ 16.0\\ 18.2}&
\tabincell{c}{ 388.37\\ 219.48\\ 101.49\\ 54.19\\ 28.59} &
\tabincell{c}{ 27\\ 35\\ 39 \\ 47\\ 55} &
\tabincell{c}{ 5.70\\ 5.22\\ 4.49\\ 3.91\\ 3.39} \\
\hline
\end{tabular}
\caption{Strong scalability study on the RANS2 problem using parallel graph partitioning.}\label{tab:cfd_zoltan}
\end{center}
\end{table}

\begin{table}[!h]
\begin{center}
\begin{footnotesize}
\begin{tabular}{c|c|c|c|c|c|c}
Matrix  & Method   & \begin{tabular}{c}Factorization \\ time (s)\end{tabular} & \begin{tabular}{c}Solving \\ time (s)\end{tabular} & \begin{tabular}{c}Total \\ time (s)\end{tabular} & Its &  Mem \\
\hline
RANS3 &
\tabincell{c}{pARMS\\BJ+VBARMS\\ BJ+VBILUT}&
\tabincell{c}{-\\ 99.41 \\ 20.45} &
\tabincell{c}{-\\ 187.95\\ 8997.82} &
\tabincell{c}{-\\ 287.36 \\ 9018.27} &
\tabincell{c}{-\\ 154\\ 979} &
\tabincell{c}{6.63 \\ 4.40 \\ 13.81}\\
\hline
RANS4 &
\tabincell{c}{pARMS\\BJ+VBARMS\\ BJ+VBILUT}&
\tabincell{c}{-\\ 12.05 \\ 1.16 } &
\tabincell{c}{-\\ 105.89\\ 295.20} &
\tabincell{c}{-\\ 117.94\\ 296.35} &
\tabincell{c}{-\\ 223\\ 472} &
\tabincell{c}{5.38 \\ 3.91\\ 5.26}\\
\hline
\end{tabular}
\end{footnotesize}
\end{center}
\caption{Experiments on the DPW3 Wing-1 problem. The RANS3 test case is solved on 32 processors and the RANS4 problem on 128 processors. The dash symbol $-$ in the table means that in the GMRES iteration the residual norm is very large and the program is aborted.}\label{tab:cfd_arms}
\end{table}

\section{Conclusions}

We have presented a parallel MPI-based implementation of a new variable block multilevel ILU factorization preconditioner for solving general nonsymmetric linear systems. One nice feature of the proposed solver is that it detects automatically exact or approximate dense structures in the coefficient matrix. It exploits this information to maximize computational efficiency. We have also introduced a modified compression algorithm that can find these approximate dense blocks structures, and requires only one simple to use parameter. The results show that the solver has nice parallel performance, also thanks to the use of a parallel graph partitioner, and it may be noticeably more robust than other state-of-the-art methods that do not exploit the fine-level block structure of the underlying matrix. 

\section{Acknowledgements}

The work of M. Sosonkina was supported in part by the Air Force Office of Scientific Research under the AFOSR award FA9550-12-1-0476, and by the National Science Foundation grants NSF/OCI—0941434, 0904782, 1047772. The authors acknowledge the Texas Advanced Computing Center (TACC) at the University of Texas at Austin for providing HPC resources that have contributed to the research results reported in this paper. URL: {\url{http://www.tacc.utexas.edu}}. The authors are grateful to Sven Baars for his assistance in implementing some of the algorithms described in the paper, and to the reviewers for their insightful comments that helped much improve the presentation.


\begin{thebibliography}{10}

\bibitem{ashc:95}
C.~Ashcraft.
\newblock Compressed graphs and the minimum degree algorithm.
\newblock {\em SIAM J. Scientific Computing}, 16(6):1404--1411, 1995.

\bibitem{OAxelsson_PSVassilevski_1989a}
O.~Axelsson and P.~S. Vassilevski.
\newblock Algebraic multilevel preconditioning methods, {I}.
\newblock {\em Numer. Math.}, 56:157--177, 1989.

\bibitem{OAxelsson_PSVassilevski_1990a}
O.~Axelsson and P.~S. Vassilevski.
\newblock Algebraic multilevel preconditioning methods. {II}.
\newblock {\em SIAM J. Numer. Anal.}, 27:1569--1590, 1990.

\bibitem{ZoltanHomePage}
Erik Boman, Karen Devine, Lee~Ann Fisk, Robert Heaphy, Bruce Hendrickson, Vitus
  Leung, Courtenay Vaughan, Umit Catalyurek, Doruk Bozdag, and William
  Mitchell.
\newblock {Zoltan} home page.
\newblock \url{http://www.cs.sandia.gov/Zoltan}, 1999.

\bibitem{ijcfd:aldo}
A.~Bonfiglioli.
\newblock Fluctuation splitting schemes for the compressible and incompressible
  {E}uler and {N}avier-{S}tokes equations.
\newblock {\em IJCFD}, 14:21--39, 2000.

\bibitem{NGILU}
E.F.F. Botta, A.~van~der Ploeg, and F.W. Wubs.
\newblock Nested grids {ILU}-decomposition {(NGILU)}.
\newblock {\em Journal of Computational and Applied Mathematics}, 66:515--526,
  1996.

\bibitem{ppam2014}
B.~Carpentieri, J.~Liao, and M.~Sosonkina.
\newblock {\em Parallel Processing and Applied Mathematics}, volume 8385 of
  {\em Lecture Notes in Computer Science}, chapter Variable block multilevel
  iterative solution of general sparse linear systems, pages 520--530.
\newblock In R. Wyrzykowski, J. Dongarra, K. Karczewski, and J. Wasniewski.
  Springer-Verlag., 2014.

\bibitem{VBARMS}
B.~Carpentieri, J.~Liao, and M.~Sosonkina.
\newblock {VBARMS}: {A} variable block algebraic recursive multilevel solver
  for sparse linear systems.
\newblock {\em Journal of Computational and Applied Mathematics}, 259
  (A):164--173, 2014.

\bibitem{dddh:90}
J.J. Dongarra, J.~Du Croz, I.~S. Duff, and S.~Hammarling.
\newblock A set of level 3 basic linear algebra subprograms.
\newblock {\em ACM Trans. Math. Softw.}, 16:1--17, 1990.

\bibitem{geli:81}
A.~George and J.~W. Liu.
\newblock {\em Computer Solution of Large Sparse Positive Definite Systems}.
\newblock Prentice-Hall, Englewood Cliffs, New Jersey, 1981.

\bibitem{guge:10}
A.~Gupta and T.~George.
\newblock Adaptive techniques for improving the performance of incomplete
  factorization preconditioning.
\newblock {\em SIAM J. Sci. Comput.}, 32(1):84--110, 2010.

\bibitem{metis4}
G.~Karypis and V.~Kumar.
\newblock Metis: A software package for partitioning unstructured graphs,
  partitioning meshes, and computing fill-reducing orderings of sparse matrices
  version 4.0.
\newblock \url{http://glaros.dtc.umn.edu/gkhome/views/metis}.
\newblock University of Minnesota, Department of Computer Science / Army HPC
  Research Center Minneapolis, MN 55455.

\bibitem{itsol:09}
Na~Li, B.~Suchomel, D.~Osei-Kuffuor, and Y.~Saad.
\newblock {ITSOL:} iterative solvers package.

\bibitem{pARMS}
Z.~Li, Y.~Saad, and M.~Sosonkina.
\newblock p{ARMS}: a parallel version of the algebraic recursive multilevel
  solver.
\newblock {\em Numerical Linear Algebra with Applications}, 10:485--509, 2003.

\bibitem{nesz:90}
J.~O'Neil and D.B. Szyld.
\newblock A block ordering method for sparse matrices.
\newblock {\em SIAM J. Scientific and Statistical Computing}, 11(5):811--823,
  1990.

\bibitem{quva:99}
A.~Quarteroni and A.~Valli.
\newblock {\em Domain decomposition methods for partial differential
  equations.}
\newblock Clarendon Press Oxford, 1999.

\bibitem{saad:93}
Y.~Saad.
\newblock A flexible inner-outer preconditioned {GMRES} algorithm.
\newblock {\em SIAM J. Scientific and Statistical Computing}, 14:461--469,
  1993.

\bibitem{Saad:1996:IME}
Y.~Saad.
\newblock {ILUM}: {A} multi-elimination {ILU} preconditioner for general sparse
  matrices.
\newblock {\em SIAM J. Scientific Computing}, 17(4):830--847, 1996.

\bibitem{VBILU}
Y.~Saad.
\newblock Finding exact and approximate block structures for ilu
  preconditioning.
\newblock {\em SIAM J. Sci. Comput.}, 24(4):1107--1123, 2002.

\bibitem{SAAD-SHORT}
Y.~Saad.
\newblock {\em Iterative Methods for Sparse Linear Systems.}
\newblock SIAM, 2nd edition, 2003.

\bibitem{ARMS}
Y.~Saad and B.~Suchomel.
\newblock {ARMS}: An algebraic recursive multilevel solver for general sparse
  linear systems.
\newblock {\em Numerical Linear Algebra with Applications}, 9(5):359--378,
  2002.

\bibitem{spal:94}
P.R. Spalart and S.R. Allmaras.
\newblock {A one-equation turbulence model for aerodynamic flows}.
\newblock {\em La Recherche-Aerospatiale}, 1:5--21, 1994.

\bibitem{IMF}
N.~Vannieuwenhoven and K.~Meerbergen.
\newblock {IMF}: {A}n incomplete multifrontal {LU}-factorization for
  element-structured sparse linear systems.
\newblock {\em SIAM J. Sci. Comput.}, 35(1):A270--A293, 2013.

\bibitem{vuduc:2007}
S.~Williams, L.~Oliker, R.~Vuduc, J.~Shalf, K.~Yelick, and J.~Demmel.
\newblock Optimization of sparse matrix-vector multiplication on emerging
  multicore platforms.
\newblock In {\em Proc.~ACM/IEEE Conf. Supercomputing (SC)}, 2007.

\bibitem{wozi:08}
P.~Wong and D.~Zingg.
\newblock Three-dimensional aerodynamic computations on. unstructured grids
  using a newton-krylov approach.
\newblock {\em Computers \& Fluids}, 37:107--120, 2008.

\bibitem{drag}
Drag~Prediction Workshop.
\newblock
  URL:\url{http://aaac.larc.nasa.gov/tsab/cfdlarc/aiaa-dpw/Workshop3/workshop3.html}.

\end{thebibliography}
\end{document}